# Estimation of Mean in Presence of Non Response Using Exponential Estimator

Rajesh Singh, Mukesh Kumar, Manoj K. Chaudhary

Department of Statistics, B.H.U., Varanasi (U.P.)-India

(rsinghstat@yahoo.com)

Florentin Smarandache
Chair of Department of Mathematics, University of New Mexico, Gallup, USA

(smarand@unm.edu)

Abstract

This paper considers the problem of estimating the population mean using information on auxiliary variable in presence of non response. Exponential ratio and exponential product type estimators have been suggested and their properties are studied. An empirical study is carried out to support the theoretical results.

## 1. Introduction

In surveys covering human populations, information is in most cases not obtained from all the units in the survey even after some call-backs. Hansen and Hurwitz(1946) considered the problem of non response while estimating the population mean by taking a sub sample from the non respondent group with the help of some extra efforts and an estimator was proposed by combining the information available from response and non response groups. In estimating population parameters like the mean, total or ratio, sample survey experts sometimes use auxiliary information to improve precision of the estimates. When the population mean $\overline{X}$ of the auxiliary variable x is known and in presence of non response, the problem of estimation of population mean $\overline{Y}$ of the study variable y has been discussed by Cochran (1977), Rao (1986), Khare and Srivastava(1997) and Singh and Kumar (2008). In Hansen and Hurwitz(1946) method, questionnaires are mailed to all the respondents included in a sample and a list of non respondents is prepared after the deadline is over. Then a sub sample is drawn from the set of non respondents and a direct interview is conducted with the selected respondents and the necessary information is collected.

Assume that the population is divided into two groups, those who will not respond called non–response class. Let $N_1$ and $N_2$ be the number of units in the population that belong to the response class and the non-response class respectively ($N_1+N_2=N$). Let $n_1$ be the number of units responding in a simple random sample of size n drawn from the population and $n_2$ the number of units not responding in the sample. We may regard the sample of $n_1$ respondents as a simple random sample from the response class and the sample of $n_2$ as a simple random sample from the non-response class. Let $h_2$ denote the size of the sub sample from $n_2$ non-respondents to be interviewed and $f = \frac{n_2}{h_2}$. Let $\bar{y}_1$ and $\bar{y}_{h2}$ denote the sample means of y character based on $n_1$ and $h_2$ units respectively. The estimator proposed by Hansen and Hurwitz (1946) is given by-

$$\bar{y}^* = \frac{n_1\bar{y}_1 + n_2\bar{y}_{h2}}{n} \qquad (1.1)$$

The estimator $\bar{y}^*$ is unbiased and has variance

$$V(\bar{y}^*) = (\frac{1}{n} - \frac{1}{N})S_y^2 + (f-1)\frac{N2}{N}.\frac{S_{y2}^2}{n} \qquad (1.2)$$

The population mean square of the character y is denoted by $S_y^2$ and the population mean square of y for $N_2$ non-response units of the population is denoted by $S_{y2}^2$.

Bahl and Tuteja (1991) introduced an exponential ratio –type estimator for population mean as given by

$$\bar{y}_{er} = \bar{y}\exp\left[\frac{\bar{X}-\bar{x}}{\bar{X}+\bar{x}}\right] \tag{1.3}$$

and exponential product- type estimator as

$$\bar{y}_{ep} = \bar{y}\exp\left[\frac{\bar{X}-\bar{x}}{\bar{X}+\bar{x}}\right] \tag{1.4}$$

The objective of this paper is to study Bahl and Tuteja (1991) exponential ratio- type and product- type estimators in presence of non-response.

**2. Suggested estimator**

First we assume that the non response is only on study variable. The estimator $\bar{y}_{er}$ under non response will take the form

$$\bar{y}^*_{er} = \bar{y}^* \exp\left[\frac{\bar{X}-\bar{x}}{\bar{X}+\bar{x}}\right] \tag{2.1}$$

To obtain the bias and MSE of the estimator $\bar{y}^*_{er}$ we write

$\bar{y}^* = \bar{Y}(1+e^*{}_0), \quad \bar{x} = \bar{X}(1+e_1)$

Such that

$E(e_0^*)=E(e_1)=0,$

and

$$\left.\begin{aligned} E(e_0^2) &= \frac{V(\bar{y}^*)}{\bar{Y}^2} \\ E(e_1^2) &= \frac{V(\bar{x})}{\bar{X}^2} = (\frac{1}{n}-\frac{1}{N})\frac{S_x^2}{\bar{X}^2} \\ E(e_o^* e_1) &= \frac{\text{cov}(\bar{y}^*,\bar{x})}{\bar{Y}\bar{X}} = (\frac{1}{n}-\frac{1}{N})\frac{S_{XY}}{\bar{X}\bar{Y}} \end{aligned}\right\} \quad (2.2)$$

where $S_{xy} = \frac{1}{(N-1)}\sum_{j=1}^{N}(x_j - \bar{X})(y_j - \bar{Y})$.

Now expressing $\bar{y}_{er}^*$ in terms of e's we have

$$\bar{y}_{er}^* = \bar{Y}(1+e_0)\exp\left\{\frac{-e_1}{2(1+\frac{e_1}{2})}\right\} \quad (2.3)$$

Expanding the right hand side of (2.3) and neglecting the terms involving powers of e's greater than two we have

$$\bar{y}_{er}^* = \bar{Y}(1+e_0-\frac{e_1}{2}-\frac{e_0e_1}{2}+\frac{3e_1^2}{8}) \quad (2.4)$$

Taking expectations of both sides of (2.4), we get the bias of $\bar{y}_{er}^*$ to the first degree of approximation, as

$$B(\bar{y}_{er}^*) = (\frac{1}{n}-\frac{1}{N})\bar{Y}\left[\frac{3C_x^2}{8}-\frac{\rho C_y C_x}{2}\right] \quad (2.5)$$

Squaring both sides of (2.4) and neglecting terms of e's involving powers greater than two we have

$$(\bar{y}_{er}^{*}-\bar{Y})^2=\bar{Y}^2\{e_0^2+\frac{e_1^2}{4}-e_0e_1\} \qquad (2.6)$$

Taking expectations of both sides of (2.6) we get the MSE (to the first degree of approximation) as

$$MSE(\bar{y}_{er}^{*})=(\frac{1}{n}-\frac{1}{N})\bar{Y}^2[C_y^2+\frac{C_x^2}{4}-\rho C_y C_x]+\frac{(f-1)N_2}{nN}S_{y2}^2 \qquad (2.7)$$

**3. Exponential product type estimator**

The estimator $\bar{y}_{ep}$ under non response (only on study variable) will take the form

$$\bar{y}_{ep}^{*}=\bar{y}^{*}\exp\left[\frac{\bar{x}-\bar{X}}{\bar{x}+\bar{X}}\right] \qquad (3.1)$$

Following the procedure of section 2, we get the bias and MSE of $\bar{y}_{ep}^{*}$ as

$$Bias(\bar{y}_{ep}^{*})=(\frac{1}{n}-\frac{1}{N})\bar{Y}\left[\frac{-C_x^2}{8}+\frac{\rho C_x C_y}{2}\right] \qquad (3.2)$$

and

$$MSE(\bar{y}_{er}^{*})=(\frac{1}{n}-\frac{1}{N})\bar{Y}^2[C_y^2+\frac{C_x^2}{4}+\rho C_y C_x]+\frac{(f-1)N_2}{nN}S_{y2}^2 \qquad (3.3)$$

**4. Non response on both y and x**

We assume that the non response is both on study and auxiliary variable. The estimator $\bar{y}_{er}$ and $\bar{y}_{ep}$ under non response on both the variables takes the following form respectively-

$$\bar{y}_{er}^{**}=\bar{y}^{*}\exp\left[\frac{\bar{X}-\bar{x}^{*}}{\bar{X}+\bar{x}^{*}}\right] \qquad (4.1)$$

$$\bar{y}_{ep}^{**} = \bar{y}^{*} \exp\left[\frac{\bar{x}^{*} + \bar{X}}{\bar{x}^{*} + \bar{X}}\right] \tag{4.2}$$

To obtain the bias and MSE of the estimator $\bar{y}_{er}^{**}$ and $\bar{y}_{ep}^{**}$ we write

$\bar{y}^{*} = \bar{Y}(1+e_0)$ , $\bar{x}^{*} = \bar{X}(1+e_1^{*})$ , $\bar{x} = \bar{X}(1+e_1)$

Such that

$E(e_0)=E(e_1^{*})=0$

$$\left.\begin{aligned} E(e_1^{*2}) &= \frac{v(\bar{x}^{*})}{\bar{X}^2} = (\frac{1}{n}-\frac{1}{N})S_x^2 + (f-1)\frac{N_2}{N}.\frac{S_{x2}^2}{n} \\ E(e_0 e_2) &= \frac{cov(\bar{y}^{*},\bar{x}^{*})}{\overline{YX}} \\ &= \frac{1}{\overline{YX}}\left[(\frac{1}{n}-\frac{1}{N})S_{xy} + (f-1)\frac{N_2}{N}.\frac{S_{xy}}{n}\right] \end{aligned}\right\} \tag{4.3}$$

The population mean square of the character x is denoted by $S_x^2$ and the population mean square of x for $N_2$ non response units of the population is denoted by $S_{x2}^2$.

The biases and MSE of the estimators $\bar{y}_{er}^{**}$ and $\bar{y}_{ep}^{**}$ are given by $\bar{y}_{er}^{**}$ and $\bar{y}_{ep}^{**}$

$$B(\bar{y}_{er}^{**}) = (\frac{1}{n}-\frac{1}{N})\bar{Y}\left[\frac{3}{8}.\frac{S_x^2}{\bar{X}^2} - \frac{S_{xy}}{2\overline{YX}}\right] + \bar{Y}\frac{(f-1)}{n}.\frac{N_2}{N}\left[\frac{3}{8}.\frac{S_{x2}^2}{\bar{X}^2} - \frac{S_{xy2}}{2\overline{YX}}\right]$$

$$= (\frac{1}{n}-\frac{1}{N})\bar{Y}\left[\frac{3}{8}.C_x^2 - \frac{\rho C_x C_y}{2}\right] + \bar{Y}\frac{(f-1)}{n}.\frac{N_2}{N}\left[\frac{3}{8}.C_x'^2 - \frac{\rho_2 C_x' C_y'}{2}\right] \tag{4.4}$$

where $C_x^{'} = \frac{S_{x2}}{\overline{X}}$ , $C_y^{'} = \frac{S_{y2}}{\overline{Y}}$ .

$$B(\overline{y}_{ep}^{**}) = (\frac{1}{n} - \frac{1}{N})\overline{Y}\left[\frac{-1}{8}.\frac{S_x^2}{\overline{X}^2} + \frac{S_{xy}}{2\overline{Y}\overline{X}}\right] + \overline{Y}\frac{(f-1)}{n}.\frac{N_2}{N}\left[\frac{-1}{8}.\frac{S_{x2}^2}{\overline{X}^2} + \frac{S_{xy2}}{2\overline{Y}\overline{X}}\right]$$

$$= (\frac{1}{n} - \frac{1}{N})\overline{Y}\left[\frac{-1}{8}.C_x^2 + \frac{\rho C_x C_y}{2}\right] + \overline{Y}\frac{(f-1)}{n}.\frac{N_2}{N}\left[\frac{-1}{8}.C_x^{'2} + \frac{\rho_2 C_x^{'} C_y^{'}}{2}\right] \quad (4.5)$$

$$MSE(\overline{y}_{er}^{**}) = (\frac{1}{n} - \frac{1}{N})\overline{Y}^2\left[C_y^2 + \frac{C_x^2}{4} - \rho C_x C_y\right]$$

$$+ \quad \overline{Y}^2\frac{(f-1)}{n}.\frac{N_2}{N}\left[C_y^{'2} + \frac{C_x^{'}}{4} - \rho_2 C_x^{'} C_y^{'}\right] \quad (4.6)$$

$$MSE(\overline{y}_{ep}^{**}) = (\frac{1}{n} - \frac{1}{N})\overline{Y}^2\left[C_y^2 + \frac{C_x^2}{4} + \rho C_x C_y\right]$$

$$+ \quad \overline{Y}^2\frac{(f-1)}{n}.\frac{N_2}{N}\left[C_y^{'2} + \frac{C_x^{'}}{4} + \rho_2 C_x^{'} C_y^{'}\right] \quad (4.7)$$

From expressions (2.7),(3.3),(4.6),(4.7), we observe that the MSE expressions of suggested estimators have an additional term (which depends on non-response)as compaired to the estimator proposed by Bahl and Tuteja ((1991)(without non response).

**5. Efficiency comparisons:**

From (1.2), (2.7), (3.3), (4.6) and (4.7), we have

First we compare the efficiencies of $\overline{y}_{er}^{*}$ and $\overline{y}^{*}$

$$MSE(\overline{y}_{er}^{*}) - V(\overline{y}^{*}) \le 0$$

$$\frac{C_x^2}{4} - \rho C_x C_y \le 0$$

$$\rho \ge \frac{C_x}{4C_y} \qquad (5.1)$$

When this condition is satisfied $\bar{y}_{er}^{*}$ will be better estimator than $\bar{y}^{*}$.

Next we compare the efficiencies of $\bar{y}_{ep}^{*}$ and $\bar{y}^{*}$

$$MSE(\bar{y}_{ep}^{*}) - V(\bar{y}^{*}) \le 0$$

$$\frac{C_x^2}{4} + \rho C_x C_y \le 0$$

$$\rho \le -\frac{C_x}{4C_y} \qquad (5.2)$$

When this condition is satisfied $\bar{y}_{ep}^{*}$ will be better estimator than $\bar{y}^{*}$.

Next, we compare the efficiencies of $\bar{y}_{er}^{**}$ and $\bar{y}^{*}$

$$MSE(\bar{y}_{er}^{**}) - V(\bar{y}^{*}) \le 0$$

$$(N-n)\alpha \le (1-f)N_2\alpha'$$

$$\frac{\alpha}{\alpha'} \le \frac{(1-f)N_2}{N-n} \qquad (5.3)$$

where $\alpha = \frac{C_x^2}{4} - \rho C_y C_x$, $\quad \alpha' = \frac{C_x^2}{4} - \rho C_y' C_x'$

When this condition is satisfied $\bar{y}_{er}^{**}$ will be better estimator than $\bar{y}^{*}$.

Finally we compare the efficiencies of $\bar{y}_{ep}^{**}$ and $\bar{y}^{*}$

$$MSE(\bar{y}_{ep}^{**}) - V(\bar{y}^{*}) \le 0$$

$$(N-n)\lambda \le (1-f)N_2\lambda'$$

$$\frac{\lambda}{\lambda'} \le \frac{(1-f)N_2}{N-n} \qquad (5.4)$$

where $\lambda = \frac{C_x^2}{4} + \rho C_y C_x$, $\quad \lambda' = \frac{C_x^2}{4} + \rho C_y' C_x'$

When this condition is satisfied $\bar{y}_{ep}^{**}$ will be better estimator than $\bar{y}^{*}$.

**6. Empirical study**

For numerical illustration we consider the data used by Khare and Sinha (2004,p.53). The values of the parameters related to the study variate y (the weight in kg) and the auxiliary variate x (the chest circumference in cm) have been given below.

$\bar{Y} = 19.50$ $\quad \bar{X} = 55.86$ $\quad S_y = 3.04$

$S_x = 3.2735$ $\quad S_{y2} = 2.3552$ $\quad S_{x2} = 2.51$

$\rho = 0.85$ $\quad \rho_2 = 0.7290$

$N_1 = 71$ $\quad N_2 = 24$ $\quad N = 95$ $\quad n = 35$

Here, we have computed the percent relative efficiencies (PRE) of different suggested estimators with respect to usual unbiased estimator $\bar{y}^{*}$ for different values of f.

**Table 5.1: PRE of different proposed estimators**

| Values of w | f values | $\bar{y}^*$ | $\bar{y}^*_{er}$ | $\bar{y}^*_{ep}$ | $\bar{y}^{**}_{er}$ | $\bar{y}^{**}_{ep}$ |
|---|---|---|---|---|---|---|
| 0.10 | 1.50 | 100 | 263.64 | 45.47 | 263.65 | 45.47 |
| | 2.00 | 100 | 263.62 | 45.47 | 263.65 | 45.47 |
| | 2.50 | 100 | 263.61 | 45.48 | 263.65 | 45.47 |
| | 3.00 | 100 | 263.59 | 45.48 | 263.64 | 45.47 |
| 0.20 | 1.50 | 100 | 263.62 | 45.47 | 263.65 | 45.47 |
| | 2.00 | 100 | 263.59 | 45.48 | 263.64 | 45.47 |
| | 2.50 | 100 | 263.56 | 45.48 | 263.64 | 45.47 |
| | 3.00 | 100 | 263.52 | 45.48 | 263.63 | 45.47 |
| 0.30 | 1.50 | 100 | 263.61 | 45.48 | 263.65 | 45.47 |
| | 2.00 | 100 | 263.56 | 45.48 | 263.64 | 45.47 |
| | 2.50 | 100 | 263.51 | 45.48 | 263.63 | 45.47 |
| | 3.00 | 100 | 263.46 | 45.48 | 263.62 | 45.47 |

From table 5.1, we conclude that the estimators which use auxiliary information performs better than Hansen and Hurwitz(1946) estimators $\bar{y}^*$. Also when non response rate increases, the efficiencies of suggested estimators decreases.